\documentclass{primus}

\usepackage{amsmath,amssymb}
\usepackage{array}
\usepackage{color}
\usepackage{graphicx}
\usepackage{url}
\usepackage{longtable}
\usepackage{rotating}
\usepackage{verbatim}
\newcolumntype{x}[1]{>{\raggedright\let\newline\\\arraybackslash\hspace{0pt}}p{#1}}
\newcolumntype{y}[1]{>{\centering\let\newline\\\arraybackslash\hspace{0pt}}p{#1}}

\makeatletter
\g@addto@macro{\UrlBreaks}{\UrlOrds}
\makeatother

\title{A ``Rule-of-Five'' Framework for Models and Modeling to Unify Mathematicians and Biologists and Improve Student Learning}

\author{C.~D.~Eaton\\
Environmental Literacy Program\\
Unity College\\
Unity, ME 04444, USA\\
ceaton@unity.edu
\and
H.~C.~Highlander\\
Department of Mathematics\\
University of Portland\\
Portland, OR 97203, USA\\
highland@up.edu
\and
K.~D.~Dahlquist\\
Department of Biology\\
Loyola Marymount University\\
Los Angeles, CA 90045, USA\\
kdahlquist@lmu.edu
\and
M.~D.~LaMar\\
Department of Biology\\
College of William and Mary\\
Williamsburg, VA  23185, USA\\
mdlama@wm.edu
\and
G.~Ledder\\
Department of Mathematics\\
University of Nebraska-Lincoln\\
Lincoln, NE 68588-0130, USA\\
gledder@unl.edu
\and
R.~C.~Schugart\\
Department of Mathematics\\
Western Kentucky University\\
Bowling Green, KY 42101, USA\\
richard.schugart@wku.edu
}


\newcommand{\AmSLaTeX}{$\cal A$\kern-.1667em\lower.5ex\hbox{$\cal
M$}\kern-.125em $\cal S$-\LaTeX}

\begin{document}


\makePtitlepage
\makePtitle

\begin{abstract}

Despite widespread calls for the incorporation of mathematical modeling into the undergraduate biology curriculum, there is lack of a common understanding around the definition of modeling, which inhibits progress. In this paper, we extend the ``Rule of Four,'' initially used in calculus reform efforts, to a framework for models and modeling that is inclusive of varying disciplinary definitions of each. This unifying framework allows us to both build on strengths that each discipline and its students bring, but also identify gaps in modeling activities practiced by each discipline. We also discuss benefits to student learning and interdisciplinary collaboration.

\end{abstract}

\listkeywords

%
%
%
%
%
%
%

\section{MANY CALLS TO ACTION}
From computer games and medicine to weather predictions and new technologies, nearly every aspect of our lives is influenced by mathematical modeling.  Primary barriers to forward progress in teaching mathematical modeling across our partner disciplines are misconceptions and biases around what constitutes modeling.  We (the co-authors) are part of an interdisciplinary working group at the National Institute for Mathematical and Biological Synthesis (NIMBioS) that has brought together mathematicians, biologists, and education researchers to address teaching quantitative biology, especially modeling.  Each of us has experienced reaching out to collaborate with a member of another discipline, only to have the conversation shut down before it has really started because of the naive assumptions we make about the other's discipline.  We posit that if mathematicians and biologists alike can improve their understanding of the similarities and differences in their approaches to and language around modeling, then each discipline will play a more effective role in advancing the other \cite{Sorgo10}, and we will be able to teach this valuable skill more effectively.  In this paper we describe a framework for models and modeling that can bridge the communication gap between disciplinary boundaries, enabling mathematicians, statisticians, and biologists to come together to improve student learning.

In 2012, the President's Council of Advisors on Science and Technology (PCAST) released the \textit{Engage to Excel} report \cite{PCAST12}. The report recommends that we should engage in a national experiment that encourages faculty from math-intensive fields other than math to be involved in teaching mathematics as a way to help close the mathematics achievement gap \cite{PCAST12}. The subtext is that traditionally trained mathematics educators are failing at helping our students succeed at mathematics as applied to science and technology. Soon after, the National Research Council (NRC) released a study, titled \textit{The Mathematical Sciences in 2025}, which suggests that in order for the mathematical sciences to remain strong in the United States, the education of students should be conducted in a cross-disciplinary manner that reflects these ever-changing realities \cite{NRC13}. This requires a rethinking of the curricula in the mathematical sciences, especially for mathematics and statistics departments, in order to provide the additional quantitative skills needed for students entering the workforce in the fastest growing career fields, such as those in STEM.

Mathematics professional societies responded to this call with the \textit{Common Vision} project, which identifies ways of improving undergraduate curricula and education in the mathematical sciences by bringing together leaders from five mathematical and statistical associations \cite{CommonVision}.  The project summarizes the collective recommendations of seven other curricular guides on undergraduate education from the associations, as well as adds its own suggestions for improving undergraduate education in the mathematical sciences, especially in the first two years of college.  The \textit{Common Vision} project concludes that departments should increase efforts to update curricula, support evidence-based pedagogical methods, and establish connections with other disciplines.

From the survey of curricular guides, the \textit{Common Vision} project identifies six themes for improving undergraduate curricula.  They are:  (1) to find more pathways into and through the curriculum for both STEM and non-STEM majors; (2) to increase the presence of statistics in student training; (3) to increase the use of modeling and computation in order to enhance conceptual understanding and introduce the scientific method into math classes; (4) to connect to skills needed in other disciplines; (5) to improve communication skills through technical writing and presentations; and (6) to aid in the transition from secondary to post-secondary education as well as from two-year to four-year institutions for transfer students.

Modeling can play an important part of several of these themes, not just where it is mentioned explicitly.  \textit{Common Vision} notes that an early introduction to modeling, along with statistics and computation, can be a pathway ``into and through mathematical sciences curricula \cite{CommonVision}.''  The Society for Industrial and Applied Mathematics (SIAM) suggests that professional societies should play a greater role in the incorporation of modeling throughout the undergraduate curriculum \cite{SIAM14}. This is reflected in the newly formed SIAM Special Interest Activity Group on Applied Mathematics Education, which recommends the development of a first-year modeling course that ``precedes and motivates the study of calculus and other fundamental mathematics for STEM majors \cite{SIAM12}.''  The American Statistical Association's Curriculum Guidelines for Undergraduate Programs in Statistical Science suggests that incorporating statistical modeling with simulations into mathematics courses can improve computational skills \cite{SIAM14}.  The Mathematics Association of America's subcommittee on Curriculum Renewal Across the First Two Years (CRAFTY) released a report in 2004 \cite{CRAFTY}, titled \textit{The Curriculum Foundations Project: Voices of the Partner Disciplines}, which emphasizes mathematical modeling in math courses.  The report notes that ``every disciplinary group in every workshop'' identifies mathematical modeling as an essential part of training students in the first two years of their undergraduate experience.  Furthermore, having students engage in mathematical modeling can ``provide a mechanism for communication, expression, and reasoning that is cross-cultural and cross-disciplinary \cite{CRAFTY},'' while having students develop a ``set of transferable skills that has the potential to be far more impactful on their futures \cite{Garfunkel16}.''

These findings are not exclusive to mathematics communities, but in fact run parallel to findings within mathematics' partner disciplines. For example, in biology, three reports, \textit{Bio2010} by the National Academy of Science (NAS), \textit{Vision and Change} by the American Association for the Advancement in Science (AAAS) and the \textit{Scientific Foundations for Future Physicians (SFFP)} Report of the American Association of Medical Colleges - Howard Hughes Medical Institute Committee, all mention the important role of mathematics, and specifically modeling, in the future of biology as a discipline \cite{Visionandchange, SFFP, Bio2010}. \textit{Bio2010} outlines, in an incredible amount of detail, the core concepts that future research biologists need from mathematics and computing, which include multiple mentions of modeling (both mathematical and statistical) throughout, as well as a section devoted to important modeling concepts \cite{Bio2010}. \textit{Vision and Change} specifically names the \textit{ability to use modeling and simulation} as a \textit{Core Competency} (emphasis added) \cite{Visionandchange}.  The \textit{SFFP} report also identifies modeling as a \textit{Core Competency} in the following way:  students should be able to ``apply quantitative knowledge and reasoning--including integration of data, modeling, computation, and analysis--and informatics tools to diagnostic and therapeutic clinical decision making \cite{SFFP}.''

\section{WANTED: A COMMON FRAMEWORK FOR MODELS AND MODELING}

At the onset of our own interdisciplinary conversation as researchers, it was clear that there were disciplinary differences by what is meant by models and modeling between the mathematicians, the biologists, the statisticians, and the STEM educators.  We generated the following questions:  Does modeling require the use of data, or even numbers?  Does a model or the process of modeling require the use of symbolic equations or formulas?  Does a schematic qualify as a model? Are you still engaged in modeling if you have not completed an entire iterative modeling process? (Hint: we will suggest resolutions to these questions in the paper that follows.)

This linguistic confusion about what constitutes models and modeling between mathematicians and the other disciplines is a barrier to interdisciplinary conversation.  It is compounded by the historical, philosophical, and physical separation of departments of mathematics and departments of statistics on large campuses.  In addition, the teaching of statistics occurs in many different departments (statistics, mathematics, biology, psychology, economics, business, education, kinesiology, etc.) on many campuses, both large and small. Furthermore, individual biologists identify primarily with one of the many different subdisciplines of biology, each with their own approaches and rich modeling traditions, \textit{e.g.}, physiological modeling, ecological modeling, and more recently, systems biology modeling.  Lastly, there is also the specific field of mathematical biology, whose practitioners are asked to move fluidly between the identities of mathematician and biologist, while still respecting the disciplinary cultures of each.  If, in all of this diversity, we are not clear about our definitions of models and modeling, then our students will not be clear.  We suggest, however, that it is possible to articulate an overarching framework for models and modeling that will unify what seem like disparate traditions with the advantage of improving student learning.

\begin{figure}[htbp]
\centering
\includegraphics[scale=0.5]{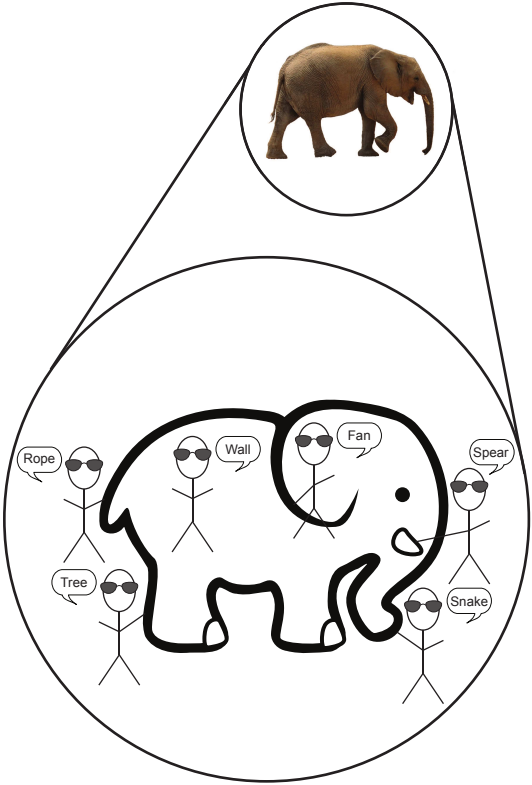}
\caption{Depiction of the blind men and the elephant parable, which illustrates the relationship between reality (the small circle) and different types of model representations (the large circle) \cite{Saxe1873}, discussed fully in section 3.1. Elephant images used in this illustration are Creative Commons public domain from Pixabay.com.}
\label{fig:elephant}
\end{figure}
\begin{figure}[htbp]
\centering
\includegraphics[scale=0.5]{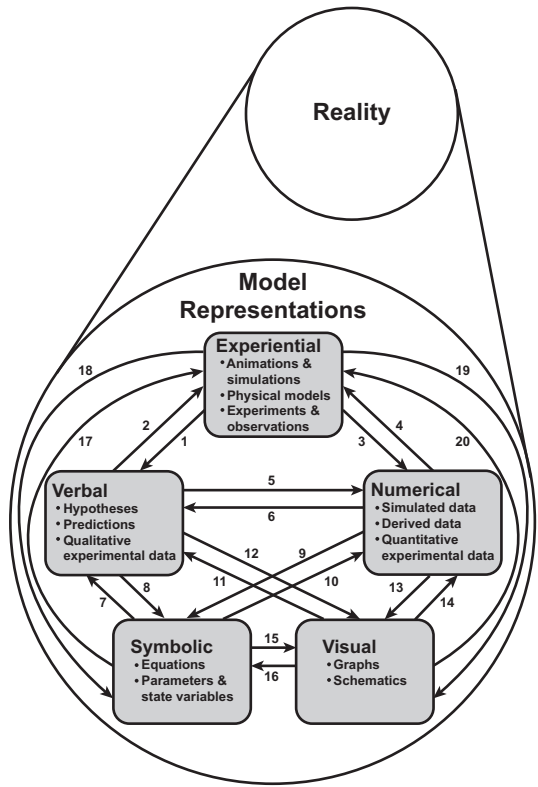}
\caption{``Rule-of-five'' framework for models and modeling. Each box is a model representation (defined in Table \ref{table:boxes}). Each arrow is an activity in the modeling process (defined in Table \ref{table:arrows}).}
\label{fig:framework}
\end{figure}

For our discussions below, we define \textit{model} as a simplified, abstract or concrete representation of relationships and/or processes in the real world, constructed for some purpose.  By \textit{simplified}, we mean that the model corresponds to a caricature of the real world rather than the real world itself, as shown in Figure \ref{fig:elephant}, which depicts the parable of the blind men and the elephant, discussed more fully in section 3.1.  The \textit{purpose} of a model is typically to enhance understanding of the process or relationship being modeled; there is a rich literature on model utility to which we refer the reader \cite{Epstein,Odenbaugh,Svoboda13}. Next, we categorize model representations into five types in our framework:  \textbf{Experiential}, \textbf{Numerical}, \textbf{Symbolic}, \textbf{Verbal}, and \textbf{Visual} (see the boxes in Figure \ref{fig:framework}).  This ``rule-of-five'' categorization has been used previously in different contexts and will be discussed in detail in Section 3 \cite{Hughes-Hallett98,Hughes-Hallett13,Simundza06}.  The rationale for making these categories explicit is that most mathematicians look at a model primarily as a collection of formulas, which is the way of looking at a model that is least accessible to biologists, thereby serving as a deterrent to the goal of increasing the amount of modeling that occurs in biology.  We seek to bridge the gap between mathematical and biological cultures by introducing the concept of multiple representations of models.  If both mathematicians and biologists appreciate that the same model given with formulas by the mathematician can be thought of in terms of graphs, data, or experiences by the biologist, then it is much easier to achieve a common understanding that is more nuanced than the individual understandings of the members of each discipline.

With model representations as \textit{objects} clearly defined, we now define \textit{modeling} as a \textit{process} composed of the various sets of activities involved in a larger modeling enterprise.  These individual \textit{modeling activities} include: 1) moving from observations of reality to an abstracted model, either as an initial step in developing a model or as part of a model revision, 2) moving from one model representation to another representation of the same model (the arrows in the framework figure, Figure \ref{fig:framework}), or 3) comparing models to each other (\textit{e.g.}, model selection) or to reality (model validation). One or more of these modeling activities comprises \textit{modeling}. Some of the tasks involved in the modeling enterprise, such as finding equilibrium solutions or solving a quadratic equation, are \textit{mathematical activities} for which knowing the context of the model is not required.  Hence, we identify \textit{modeling activities} as comprising just those tasks that make sense only in the context of modeling, such as translating a verbal representation into a mathematical representation or comparing the predictions made by different models; these are the tasks that require moving from one model representation to another representation of the same model (the arrows in the framework figure, Figure \ref{fig:framework}).

Our definition re-frames the modeling process as having two levels of detail: a holistic level that defines the \textit{modeling enterprise} as the traditional iterative process, with a goal of creating a useful final model, and a finer more granular level that considers individual tasks that, taken together, comprise a modeling process.  Our definition includes conceptions of the modeling process as a \textit{complete} set of steps that are iterated (for example, as presented in the GAIMME report \cite{Garfunkel16}: state the problem, make simplifying assumptions, mathematize, analyze the model, refine or extend the model).  While it is important for students to experience the full modeling enterprise, in particular the iterative nature of modeling -- one's first attempt at a model is seldom adequate, -- we suggest that there is value in practicing individual modeling activities as well.  Thus, defining modeling as any set of modeling activities is both more expansive and more inclusive of what ``counts'' as modeling.  Pedagogically, our definition allows the instructor to more easily scaffold modeling into the curriculum, which is especially useful in our partner disciplines and also allows us to acknowledge to our partner disciplines that \textit{we are all modelers}.  

In the sections below, we will explain what we mean by the ``rule-of-five'' framework, with special attention to \textbf{Experiential} as a representation that is critical to and often missing from how we are currently teaching modeling to students.  We then give illustrative examples of model representations, modeling activities, and modeling pathways for the logistic growth model, the Hardy-Weinberg Equilibrium (HWE) model, and a physical model of the structure of DNA. These examples are often found in introductory biology classes, and can therefore be used by students to practice modeling activities in biology classrooms and learn modeling concepts more deeply. We show how this framework can be used to 1) provide a common language with which to engage in interdisciplinary conversations around the modeling process, inclusive of whether approached from the point of view of mathematics, statistics, or biology, and 2) provide a common framework for the teaching of modeling across disciplines.  While our focus in this paper is to provide a framework that will facilitate communication between mathematicians, statisticians, and biologists, we believe this framework to be adaptable to any discipline. 

\section{``RULE-OF-FIVE'' FRAMEWORK}
\subsection{Background and Justification}

The CRAFTY report, described above, noted that by engaging in mathematical modeling, students have an opportunity to describe their work, ``\textbf{Analytically}, \textbf{Graphically}, \textbf{Numerically}, and \textbf{Verbally} \cite{Ganter04}.''  These model descriptions are historically referred to as ``the rule of four \cite{Windham08}.''  This concept came about in the 1980s and 1990s when it was recognized that over half of students enrolled in calculus courses in the United States did not finish the course \cite{Douglas86}.  This led to a number of calculus reform projects, many of which were funded by the National Science Foundation.  Several of these projects used a ``multiple representation approach,'' as was suggested in the National Council of Teachers of Mathematics' \textit{Curriculum and Evaluation Standards for School Mathematics} \cite{Utter96, NCTM89}.  These ideas then came to the forefront with the work of Hiebert and colleagues \cite{Hiebert86,Hiebert92}.  Hiebert argued that in order to learn mathematics, a student must understand mathematics.  Such understanding occurs when there is a continually evolving and strengthening network of connections to internal mental representations of mathematical ideas and procedures.  Hiebert proposed a framework to aid in the understanding of mathematics by making connections between \textbf{Numerical},  \textbf{Graphical}, and \textbf{Symbolic} representations. The \textit{Calculus Connections Project} of Oregon State University implemented a reform calculus emphasizing the representations of \textbf{Graphical}, \textbf{Numerical}, and \textbf{Symbolic}, and the importance of switching between representations \cite{Dick94,Utter96}.  The \textit{St. Olaf Project} also described moving among the representations of \textbf{Graphical}, \textbf{Numerical}, and \textbf{Algebraic} as being crucial to learning the concepts of calculus \cite{Ostebee94}.  The \textit{Calculus Consortium} at Harvard University coined the term ``rule of three'' in their reform textbook \cite{Gleason92}.  The ``rule'' stated that equal weight should be given to describing topics \textbf{Algebraically}, \textbf{Numerically}, and \textbf{Geometrically} \cite{Hughes-Hallett95}.  Note that the original word choice for the representations parallels course descriptions at the secondary school level \cite{Hughes-Hallett16}.  In the second edition of the text \cite{Hughes-Hallett98}, this became a ``rule of four'' with a fourth equal-emphasis on \textbf{Verbal} descriptions of math problems by teachers and students.  In subsequent editions of the text \cite{Hughes-Hallett05,Hughes-Hallett13}, the descriptions of the four representations became \textbf{Graphical}, \textbf{Numerical}, \textbf{Symbolic}, and \textbf{Verbal}.  Since the inception of calculus reform in the 1990s, it is now commonplace to see a wider variety of problems and the use of multiple representations in ``traditional'' math textbooks \cite{Hughes-Hallett06}.

More recently, a fifth rule was proposed by Simundza in a laboratory course for precalculus \cite{Simundza06}.  This fifth rule, \textbf{Experiential}, is, in the words of Simundza, a ``direct sensory experience of quantitative phenomena \cite{Simundza06}.''  The importance of the \textbf{Experiential} representation is alluded to in the American Mathematical Association of Two-Year Colleges (AMATYC) standards where they propose that mathematics be taught like the sciences as a laboratory discipline \cite{AMATYC}.

Introducing students to multiple forms of representation is well-documented to improve student learning.  For example, the first principle guideline from the Universal Design for Learning (UDL) is a recommendation to provide students multiple forms of representation since there are diverse ways learners comprehend information \cite{Cast11}. UDL is a framework that addresses high variability in learners' responses to instruction by suggesting flexibility in the curriculum to meet the varied needs of the students.  Using multiple representations can allow students to make connections within and between concepts \cite{Cast11}. When students solved problems using more than one representation, student performance was better than for those learners who used a single strategy \cite{Ainsworth99,Cox95,Tabachneck94}.

Our use of \textbf{experiential} in the context of modeling is slightly different than, but still in alignment with, what is meant by ``experiential learning'' in other contexts \cite{Bransford99}.  There is an abundance of evidence that experiential learning can improve student-learning outcomes \cite{Kolb05,Kolb15,Kolb01,Simundza06}. There is also early evidence that adopting this type of framework, in particular multiple model representations and movements between them, can lead to success in subsequent quantitative courses \cite{Cast11}. Some students may enter into the learning experience more comfortable with a subset of these representations due to their own ways of knowing or disciplinary identity, but the goal is that they should know all and be able to move between them. 

The use of multiple forms of representation and the ``rule-of-five'' mirrors the parable of the blind men and the elephant as seen in Figure \ref{fig:elephant}. John Godfrey Saxe's poem version of the parable describes six blind men touching different parts of an elephant in order to ``understand'' it \cite{Saxe1873}.  Each blind man compares each different elephant part to an everyday object that is similar - a wall (side), snake (trunk), spear (tusk), tree (knee), fan (ear) or rope (tail).  Just as in the parable, it is only through the use of different forms of representation that we may hope to gain the truest understanding of a problem.  Interestingly, this analogy to the parable also works when considering the different ways in which mathematicians and biologists approach problems - each is experiencing one aspect of the problem, and through communication using a common language around modeling, we can create together a clearer picture of the world.

\subsection{Model Representations (The Boxes)}

As discussed above and shown in Figure \ref{fig:framework}
and Table \ref{table:boxes}, we extend the ``rule of five'' from its original use as an aid to calculus and precalculus instruction to a useful general description of the various types of model representations that may be used in the modeling process.

\begin{center}
\begin{longtable}{|x{2.2cm}|x{3.2cm}|x{4.2cm}|}
\caption{The five types of model representations, corresponding to Figure \ref{fig:framework} with their definitions and classroom examples.}\label{table:boxes}\\
\hline
\textbf{Type of representation} & \textbf{Definition} & \textbf{Classroom examples}\\
\hline
\endfirsthead
\textbf{Experiential} &
Direct experiences that are concrete rather than abstract such as, virtual laboratories and animations, kinesthetic experiences and manipulations of physical models, actual scientific experiments. &
\textbullet~A video of bacterial growth; beanbag biology \cite{Jungck10}; virtual laboratories 
(\textit{e.g}, SimBio \cite{Bugbox}, the BUGBOX-predator virtual laboratory \cite{Bugbox});\newline
\textbullet~Physical model of the structure of DNA;\newline
\textbullet~Experiment to measure bacterial growth in the laboratory.\\                   
\hline
\textbf{Verbal} &
Hypotheses used to design experiments, predictions, assumptions used to construct mathematical models, simple descriptions of observations, qualitative experimental data. &
\textbullet~Hypothesis: a mutation in a particular gene will reduce the rate of bacterial growth because the mutation impairs DNA replication;\newline
\textbullet~Prediction: on average global temperature will increase;\newline
\textbullet~Assumption: we assume that the population is well mixed;\newline
\textbullet~Simple descriptions of observations: the rate of increase is decreasing; we
observe far more of the blue flower type then the\\
&&purple flower type.\newline
\textbullet~Qualitative data: spiciness ratings by tasters of chili peppers.\\      
\hline
\textbf{Numerical} &
Data sets collected from model-based simulations, data calculated from other data (derived data), experimentally-collected quantitative data.&
\textbullet~Numbers of infected individuals calculated from a \textbf{Symbolic} epidemic model;\newline
\textbullet~Derived data: low density growth rate and carrying capacity calculated from plotting relative growth rate versus population for logistic growth;\newline
\textbullet~Measured population counts from experiments.\\
\hline 
\textbf{Visual} &
Graphs, schematics. &
\textbullet~A graph of relative growth rate versus population;\newline
\textbullet~A schematic of an epidemic model; stock-and-flow diagrams;\newline
\textbullet~Data visualizations (e.g., histograms, scatter plots, infographics, etc.).\\
\hline 
\textbf{Symbolic} &
Formal mathematical constructs such as formulas, equations, algorithms, parameters, and state variables. &
\textbullet~Discrete difference equation for geometric growth $x_{n+1}=\lambda x_n$ and continuous differential equation for exponential growth $\dfrac{dP}{dt}=rP$;\newline
\textbullet~If in HWE, $p$ =\\
&&frequency of one allele, $p^2$ = frequency of homozygotes for that allele.\newline
\textbullet~State variable: $P(t)$ = population at time $t$ (in years);\newline\textbullet~Equation from a linear regression;\newline
\textbullet~Equation from a probability distribution.\\
\hline 
\end{longtable}
\end{center}

\normalsize
Instructors well-versed in the ``rule of four'' should immediately find our definitions of the \textbf{Numerical}, \textbf{Symbolic}, \textbf{Verbal}, and \textbf{Visual} model representations familiar, although we have broadened their definitions to be inclusive of perspectives from mathematics, statistics, and biology. \textbf{Symbolic} representations are also referred to as \textit{Analytic} by the CRAFTY report \cite{Ganter04} and as \textit{Algebraic} by the St. Olaf Project and the Harvard Calculus Consortium \cite{Ostebee94,Gleason92}. \textbf{Visual} representations, in particular $x$-$y$ plots, are called \textit{Graphical} in the ``Rule of Four.''  We suggest \textbf{Visual} is more inclusive of the emerging field of data visualization \cite{Tufte}.

\subsubsection{Experiential representation as a link between science and mathematics}
\textbf{Experiential} representations are treated as distinct from the other representations, particularly the \textbf{Visual} representation, because the modeling enterprise is not pure mathematics, but theoretical science. Mathematics is abstract, so \textbf{Visual} representations of mathematical constructions in calculus, where the ``rule of three'' originated, are almost always graphs. Science is both abstract and concrete, with direct sensory experience playing a distinct and independent role. \textbf{Experiential} representations of a model are more easily connected with the original phenomena than any of the other representations.  In a mathematical modeling process, this may involve movie-like simulations of growing populations or bean-bag tactile manipulations, \textit{e.g.}, \cite{Jungck10}.  In biology, this would include the experiments themselves.  It may seem strange to define experiments as a model, but an experiment cannot encompass the entirety of reality (Cf. blind men and elephant parable, Figure \ref{fig:elephant}).  Referring to our definition of a model, an experiment is a concrete simplification of reality, designed with some purpose in mind, that depends on a particular experimental design formulated by a scientist and carried out with some type of apparatus that interacts with the real world, but is not the real world itself.

\subsubsection{Using experiential approaches in teaching}
An advantage for students working with \textbf{Experiential} representations is that they provide an entryway for understanding models that do not require the abstraction of \textbf{Symbolic} or \textbf{Visual} representations, the integration of detail needed to understand a \textbf{Numerical} representation, or the conceptual knowledge and reading comprehension needed to understand a \textbf{Verbal} representation. 

An example of an \textbf{Experiential} representation is the BUGBOX-predator virtual laboratory \cite{Bugbox}, which simulates an agent-based model of an experiment in which one predator is given a fixed amount of time to find and consume stationary prey at some initial population density of prey. We can describe to students the \textbf{Verbal} assumptions of the Holling type 2 functional response model \cite{Holling59a,Holling59b} and show them the resulting \textbf{Symbolic} equation, \textbf{Numerical} data, and graph visualization (\textbf{Visual}). These are helpful, but they are a poor substitute for showing students the animation and letting them observe for themselves that the predator must divide its time between searching and handling, with more searching at low prey density and more handling at high prey density. As an additional bonus, no prerequisite sophistication is required to appreciate purely sensory observation \cite{Jungck10}.

\subsubsection{Five representations of the logistic growth model}
The logistic growth model is another fitting example with which to illustrate the five different types of model representations.  The \textbf{Experiential} box includes experiments measuring the growth of bacteria.  This includes conducting an actual laboratory experiment or manipulating a simulation.  However, it also includes viewing the results of a time-lapse video or animation. The quantitative population data (from the experiment or simulation) and absolute and relative growth rates (derived data) go in the \textbf{Numerical} box.  A \textbf{Symbolic} representation is the set of equations that approximately describe this relationship, \textit{i.e.}, the equation for logistic growth. \textbf{Visual} representations of the model can be obtained by plotting the data in different ways, for example, as population versus time or as per capita population rate of change as a function of population size. \textbf{Verbal} representations of logistic growth include statements about the growth rate, \textit{e.g.}, ``the relative growth rate is a decreasing linear function of population,'' or as a statement linking the mathematical model to the biological concepts and processes at play, \textit{e.g.}, ``the absolute growth rate is proportional to the population and the remaining capacity for population growth.''

\subsection{Modeling Activities (The Arrows)}
As defined earlier, modeling activities for a specific model are those that connect different representations or connect the model to the real-world scenario. Recall above our definition of modeling activities involved within the process of modeling: 1) moving from observations of reality to an abstracted model, either as an initial step in developing a model or as part of a model revision, 2) moving from one model representation to another representation of the same model (the arrows in the framework figure, Figure 2), and 3) comparing models to each other (e.g. model selection) or to reality (model validation). To be clear, we intend for the boxes in Figure 2 to stand for different representations of models, and the arrows to stand for different modeling activities within the modeling process. In the next section, we discuss the use of the framework in the context of the modeling process and enterprise.

\begin{center}
\begin{longtable}{|y{1.5cm}|x{1.8cm}|x{2.4cm}|x{3.0cm}|}
\caption{Modeling activities that correspond to the arrows in Figure \ref{fig:framework} with a description and an example.}\label{table:arrows}\\
\hline
\textbf{Arrow Number} & \textbf {Box$\rightarrow$Box} & \textbf{Arrow Description} & \textbf{Classroom Example} \\
\hline
1. & Exp$\rightarrow$Ver &
Crafting a formal scientific hypothesis (which implies mechanism) based on observations; predictions; collecting of qualitative experimental data. &
\textbullet~Students make a hypothesis based on their observations of bacterial population growth.\newline
\textbullet~A population satisfying the assumptions of Hardy-Weinberg Equilibrium (HWE) will maintain the same allele frequencies generation after generation.\\
\hline
2.& Ver$\rightarrow$Exp & 
Designing an experiment. &
\textbullet~Based on a hypothesis about bacterial growth rates, make a prediction about growth rates under different conditions. Design an appropriate experiment to support or refute\\
&&&the prediction.\newline
\textbullet~Based on a hypothesis that a population is in HWE, plan an experiment to test for HWE.\\
\hline
3. & Exp$\rightarrow$Num &
Collecting data. &
\textbullet~Collecting quantitative data from the \textbf{Experiential} simulations/ animations.\newline
\textbullet~Sampling populations over time in a field study.\newline
\textbullet~Measuring a culture of bacteria over time through spectroscopy.\\
\hline
4. & Num$\rightarrow$Exp &
Feeding experimental data into a simulation; using experimental or simulated data to make a physical model; performing model validation against new &
\textbullet~Using coordinates derived from X-ray crystallography to build a physical model of DNA.\newline
\textbullet~Using data from a numerical simulation to create an animation of a population growing.\\
&&experiments.
&\textbullet~Running a longer experiment to test when an exponential growth model of bacteria fails to match the data.\\
\hline
5. & Ver$\rightarrow$Num &
Estimating; approximating. &
\textbullet~Finding estimates for parameters in a logistic model from the literature.\newline
\textbullet~Back-of-the-envelope calculations or reasoning that allow students to test a \textbf{Verbal} prediction.\\
\hline
6. & Num$\rightarrow$Ver &
Describing patterns and trends in the data; using data to refine hypotheses. &
\textbullet~Looking for trends in a data set that may indicate the presence of a carrying capacity for the population.\newline
\textbullet~Interpreting the results of a statistical  test for HWE.\\
\hline
7. & Sym$\rightarrow$Ver &
Interpreting/ analyzing a mathematical model. &
\textbullet~$\dfrac{dP}{dt}=rP-\dfrac{rP^2}{K}$ is the rate of change of population size, positively affected\\
&&&by a net positive intrinsic birth-death rate and negatively affected by intraspecific competition over resources.\newline
\textbullet~Interpreting $p$ and $q$ in the HWE model as frequency of allele A and frequency of allele a.\\
\hline
8. & Ver$\rightarrow$Sym &
Mathematizing. &
\textbullet~The converse of \#7, going from the explanation to the \textbf{Symbolic} representation.\newline
\textbullet~Activities traditionally associated with mathematical modeling, \textit{i.e.}, formalizing the language of the \textbf{Verbal} description by assigning parameters such as carrying capacity and intrinsic growth rate, state variables such as population size,\\
&&&and writing formulas and equations such as, $\dfrac{dP}{dt}=rP-\dfrac{rP^2}{K}$.\\
\hline
9. & Num$\rightarrow$Sym &
Statistical modeling. &
\textbullet~Fitting a logistic model to a data set.\newline
\textbullet~Testing the experimentally measured data of genotype frequencies to the HWE null statistical model.\\
\hline
10. & Sym$\rightarrow$Num &
Simulating data. &
\textbullet~Simulate data from the logistic differential equation using an ODE solver and a particular parameter set.\newline
\textbullet~Performing \textit{in silico} experiments.\\
\hline
11. & Vis$\rightarrow$Ver &
Interpreting visualizations of processes or results. &
\textbullet~Going from a ``stock-and-flow/ box-and-arrow'' schematic representation of a model to the \textbf{Verbal} description, a hypothesis, or prediction to test.\\
&&&\textbullet~Interpretation of a graph of logistic growth, identifying the lag phase (early exponential), the log phase (greatest rate of growth), and the stationary phase (saturation) of the logistic growth rate.\\ 
\hline
12. & Ver $\rightarrow$Vis &
Sketching a graph or schematic that illustrates a hypothesis or observation. &
\textbullet~Sketching a graph of a population that starts off growing exponentially, but then has a carrying capacity.\newline
\textbullet~Drawing a schematic of what processes might be involved in a logistically growing population.\\
\hline
13. &Num$\rightarrow$Vis &
Graphing, visualizing data. &
\textbullet~Traditional plot of population data versus time or rate of growth versus population size.\newline
\textbullet~Creating appropriate infographics for a ``big data'' set.\\
\hline
14. & Vis$\rightarrow$Num &
Interpolating data points from a graph, using a visual modeling program such as STELLA, Insightmaker, or Simulink. &
\textbullet~Interpolating between given data points.\newline
\textbullet~Estimating data points or parameter values, such as carrying capacity, from a  graph (\textit{a.k.a.}, ``Reverse engineering'').\\
\hline
15. & Sym$\rightarrow$Vis &
Graphing or drawing a schematic of a process described by a formal mathematical model. &
\textbullet~Traditional graphing of the logistic growth curve.\newline
\textbullet~Drawing a stock-and-flow diagram from seeing SIR model equations.\\
\hline
16. & Vis$\rightarrow$Sym &
Modeling based on qualitative features of a graph or processes laid out in a schematic. &
\textbullet~Writing the equation for logistic growth of a population based on a graph knowing the carrying capacity, the initial population size, and the time at which the population is at half the carrying\\
&&&capacity or the doubling time.\newline
\textbullet~Build the \textbf{Symbolic} representation of a model based on the schematic of processes involved, such as net exponential growth and death due to competition.\\
\hline
17. & Sym$\rightarrow$Exp &
Programming an animated simulation. &
\textbullet~Programming in NetLogo a simulation illustrating the logistic growth of bacteria in a virtual petri dish.\newline
\textbullet~Creating a manipulative beanbag biology experiment to explore algebraic relationships of the HWE equation.\\
\hline 
18. & Exp$\rightarrow$Sym &
Writing mathematical equations directly based on observations (experienced modelers may not need to &
\textbullet~Recognizing logistic growth is at played based on an experiential activity and immediately writing the resulting equations.\\
&&pass through additional \textbf{Verbal} or \textbf{Visual} boxes). &\\
\hline 
19. & Exp$\rightarrow$Vis & Drawing a schematic or cartoon based on experiential observations. &
\textbullet~Drawing a sketch of DNA from  an animation or a physical model.\newline\textbullet~Drawing a process schematic or concept map of processes at play when observing the growth of a bacterial cell culture.\\
\hline
20. & Vis$\rightarrow$Exp &
Experimenting with manipulatives, assembling a physical model that replicates a process schematic, or animating a sketch or cartoon of a process. &
\textbullet~Constructing a physical model of DNA from a picture.\newline
\textbullet~Perform an experiment to replicate graphical or schematic results.\\
\hline  
\end{longtable}
\end{center} 
\normalsize
Table \ref{table:arrows} describes the modeling activities contained in the arrows between boxes in the framework figure, Figure \ref{fig:framework}.  A fully connected graph (arrows which describe the transition between any two of the model representations) is the most inclusive of all potential uses of this framework.  Some activities (arrows) would be more commonly performed than others, likely dependent on the discipline.  In the table, we have interpreted the modeling activities in the context of logistic growth, Hardy-Weinberg Equilibrium, or the structure of DNA. This will also help us in the next section, where we use these individual representations and activities as building blocks of the modeling process.

We want to emphasize that work performed entirely within one particular box alone is not a ``modeling activity'' per se, but the province of the disciplinary context within which the work is performed. ``Application'' problems in which students are given a formula and asked to use it to compute an answer without explaining the meaning is not modeling \cite{Garfunkel16}. However, more extended problems in which models are used to create simulated data (arrow 10 in framework Figure \ref{fig:framework} going from \textbf{Symbolic} to \textbf{Numerical}), and in which these results are checked against additional data (arrow 4 in framework Figure \ref{fig:framework} going from \textbf{Numerical} to \textbf{Experiential}) are a much better illustration of modeling.

\section{USING THE FRAMEWORK TO UNIFY MODELING APPROACHES}
We are certainly not the first group to describe a framework for teaching modeling. Frameworks can be more or less rigid in specifying a particular order of steps in the modeling process \cite{Bender78,Brauer01,Carson14,Ledder05,Nandi13}. Recently, the GAIMME report discusses the [mathematical] modeling process, as well as resources for teaching modeling \cite{Garfunkel16,SIAM14}. It emphasizes that the entire iterative process of modeling is flexible, \textit{i.e.}, moving back-and-forth between the different stages of model formulation and analysis, with the report focusing its discussion in the context of mathematical modeling. These more modern discussions of the modeling framework are consistent with the ``messy'' and non-linear nature of what happens in actual expert practice \cite{Understandingsciencea}.

One powerful feature of our framework is that we can explicitly acknowledge and practice a variety of activities important to the modeling process without having to engage in the full modeling enterprise.  This allows us to scaffold and reinforce activities more easily, particularly in classes that are not explicitly modeling classes, such as partner discipline classes. Furthermore, if we refer to approaches taken by partner disciplines as different uses of a larger modeling framework, then when students engage in those approaches, we are setting them up to engage in the mathematical modeling process with more ease later.

Our framework is fully compatible with these envisionings of mathematical modeling, but it is more inclusive in the following ways:  Our framework 1) encourages deliberate and thoughtful development of individual modeling activities and skills in not only mathematics, but other partner discipline classes; 2) avoids what some might call ``disciplinary microaggressions'' by providing a framework inclusive of disciplinary-specific research approaches taken by mathematicians, statisticians, and biologists, and 3) emphasizes that the \textbf{Experiential} representation is crucial for student learning, particularly in partner disciplines.  

\subsection{Multiple Modeling Pathways}
In this section, we describe the modeling process as a pathway through the modeling framework. We have chosen just one representative pathway for each discipline to discuss in detail - a mathematical modeling example using logistic growth, a statistical modeling example using Hardy-Weinberg Equilibrium, and a third example from biology using the structure of DNA.  However, we explicitly acknowledge that multiple valid modeling pathways exist even within disciplines and encourage the reader to examine his or her own pathways through the framework in class and in research. The particular examples described below have been chosen due to their ubiquity as models taught in mathematics and biology. We now explore these as opportunities for instructors to engage more deeply using the framework presented here. 

\subsubsection{A mathematical modeling pathway - logistic growth}
As we noted previously, a modeling investigation becomes more and more valuable as the number and variety of connections made between different model representations increases. In teaching mathematical modeling activities to students, it is very helpful to begin with something \textbf{Experiential}. Ideally this would use real observations of biological phenomena, as would be done in biology laboratory courses. If this is not possible, a manipulatable representation of a model would serve as a reasonable starting point (for example, in biology courses without a lab component, or in mathematics/statistics courses).  Actual experiments and simulated experiments both belong in the \textbf{Experiential} box (framework Figure \ref{fig:framework}).  Depending on the initial experience, \textbf{Experiential} observations might generate either qualitative or quantitative data, which lead to \textbf{Verbal} and \textbf{Numerical} representations of a model, respectively.

The mathematical modeling process is a pathway through the framework that includes the symbolic box. In this example, we discuss the logistic growth model and follow the path \textbf{Experiential} $\rightarrow$ \textbf{Numerical} $\rightarrow$ \textbf{Visual} $\rightarrow$ \textbf{Verbal} $\rightarrow$ \textbf{Symbolic}.  We begin with an \textbf{Experiential} representation by performing an experiment to measure the growth of a bacterial culture over time.  If this cannot be done in class, one could find a time-lapse video of the phenomenon. This approach to developing a \textbf{Numerical} model would be to use raw population data to calculate absolute and relative growth rates. One could then obtain a \textbf{Visual} representation of the model by plotting the data in different ways (for example, population size versus time, or relative growth rate versus population size).  Then these steps lead to a \textbf{Verbal} model consisting of the assumption that the growth rate is proportional to population. An alternate path would go from \textbf{Experiential} $\rightarrow$ \textbf{Verbal} $\rightarrow$ \textbf{Symbolic} when the \textbf{Symbolic} model is obtained from the qualitative observation that more parents produce more offspring.  In either pathway, the \textbf{Symbolic} representation of the model does not come directly from qualitative data (\textbf{Verbal} box) or quantitative data (\textbf{Numerical} box), but comes rather from a \textbf{Verbal} representation obtained from the initial observations. In this example, the \textbf{Verbal} representation leads to a \textbf{Symbolic} representation in the form of a differential (or difference) equation. The pathway can then be extended from \textbf{Symbolic} $\rightarrow$ \textbf{Numerical} $\rightarrow$ \textbf{Visual} by running a simulation to produce data generated by the mathematical \textbf{Symbolic} model, which can be graphed.  Or the pathway can be extended from a \textbf{Symbolic} $\rightarrow$ \textbf{Experiential} representation by creating an animation with computer graphics. The pathway can be extended from \textbf{Symbolic} directly to \textbf{Visual} by solving the equation analytically and using the resulting formula to prepare a graph. In either case, the \textbf{Visual} representation may contribute a statement of model behavior to the \textbf{Verbal} representation.  Note that in this description of these pathways, the boxes and arrows are traversed several times and the order of the traversal is not fixed.  As the students engage in each of these activities, they should be thinking about what they are doing and the connections they are making.

\subsubsection{A statistical modeling pathway - Hardy-Weinberg Equilibrium}
The previous example illustrates a modeling pathway that leads to a \textbf{Symbolic} representation of a mechanistic mathematical model (logistic growth).  The next example illustrates a modeling pathway leading to a \textbf{Symbolic} representation of a statistical model, given in the context of Hardy-Weinberg Equilibrium (HWE), which is typically taught in introductory biology courses.  Pulling from previous experiments and observations  (\textbf{Experiential} $\rightarrow$ \textbf{Verbal}), students are introduced to the five major forces in evolution: (1) selection, (2) genetic drift, (3) mutation, (4) gene flow, and (5) non-random mating \cite{Raven16}.  To arrive at the probability model (\textbf{Symbolic}) of HWE, we first assume a null \textbf{Verbal} model that no evolutionary forces are present.  In particular, suppose we have an isolated (no gene flow), infinite population (no genetic drift) where random mating, no selection, and no mutation occur.  A simplified genetic model also assumes one locus and two alleles.  Given this \textbf{Verbal} model, the \textbf{Symbolic} probability model representation is derived (HWE) stating that allele and genotype frequencies do not change in subsequent generations.  Many times this is where instruction stops, but it is informative to mention that this probability model can then be used to statistically test for the presence of HWE, which is a necessary but not sufficient condition for lack of evolutionary forces.  This is accomplished by measuring genotype frequencies in a population of interest (\textbf{Experiential} $\rightarrow$ \textbf{Numerical}) and then performing a chi-squared test to the null probability model of HWE (\textbf{Numerical} $\rightarrow$ \textbf{Symbolic}).  Results of the test are then used as evidence for the presence or absence of evolutionary forces (\textbf{Symbolic} $\rightarrow$ \textbf{Verbal}).  It should be noted that most hypothesis testing in statistics follows this pathway:  observed data collected from experiment (\textbf{Experiential} $\rightarrow$ \textbf{Numerical}) is tested against a null probability model (\textbf{Numerical} $\rightarrow$ \textbf{Symbolic}), and the results are used as a quantification of evidence for a biological hypothesis (\textbf{Symbolic} $\rightarrow$ \textbf{Verbal}).

\subsubsection{A biological modeling pathway - the structure of DNA}
Finally, let's describe a pathway followed by biologists.  One of the most famous physical models in molecular biology is the physical model of the structure of DNA created by James Watson and Francis Crick, that inspired the more abstract schematic that appears in their \textit{Nature} paper \cite{Watson53}. To build this model, they integrated both quantitative and qualitative data from multiple sources (\textbf{Numerical} and \textbf{Verbal} $\rightarrow$ \textbf{Experiential}).  This physical model enabled a hypothesis about mechanism (\textbf{Experiential} $\rightarrow$ \textbf{Verbal}): ``It has not escaped our notice that the specific pairing we have postulated immediately suggests a possible copying mechanism for the genetic material'' \cite{Watson53}. One of the co-authors of the current paper regularly uses building a chemically-correct physical model of DNA as an activity in her courses.  Students gain a much deeper understanding of the chemical structure of DNA from the direct manipulation of the physical model than from \textbf{Verbal} descriptions or \textbf{Visual} aids alone.  Note that in this example, not all of the boxes or arrows were traversed, yet it is still an example of a modeling activity that enhances student learning.

\subsection{Modeling Pathways and the Modeling Enterprise}
Ultimately, the results of a modeling investigation need to be checked against original or new qualitative and/or quantitative data (or even other mathematical models). Thus, recall the first and third parts of our definition of modeling activities which include model abstraction and revision, model validation, and model selection. These are not represented in our framework Figure \ref{fig:framework} as arrows, but are important modeling activities when engaging in the full modeling enterprise. A careful critique can almost always identify features of the biological system that are missing from the model. In biology, a large amount of stochasticity is often superimposed on deterministic phenomena, so we cannot expect a model to exactly reproduce experimental data. Our first example of a modeling pathway, using \textbf{Experiential} activities to examine logistic growth, would almost surely lead to a data set that does not precisely match the deterministic logistic growth example. Thus, an important part of the modeling process includes recognizing when a model may need to be further refined to address the question at hand.  In the course of a modeling investigation with students, the instructor's role often is to remind students to pause, validate the model, and if needed, reexamine the model assumptions or mathematization to refine it.

In the instructional setting, we can use the analogy of the blind men and the elephant (Figure \ref{fig:elephant}) to remind students that results obtained from the study of a particular mathematical model pertain only to that model.  Whether they are useful in understanding the biological setting depends on comparing and contrasting the model formulation with the corresponding biological process, and model-generated data with observed or experimental data. This is a key ingredient in encouraging students to engage in any modeling activity that is part of the modeling process. Even if students are only practicing individual modeling activities (smaller pieces of a larger modeling process), it is imperative to remind them that they are working in a conceptualized model of reality \cite{Garfunkel16,Schwonke09}. They are using caricatures of reality, not dealing with reality itself, and the assumptions and results should be critically analyzed in that light.

In our example model of logistic growth, there are a number of models that one can use for limited growth. If the per capita growth rate is constant, the population is growing exponentially $\left(\mathrm{\textit{i.e.}}, \frac{1}{P}\frac{dP}{dt}=r\right)$. However there are a number of ways in which a population can experience a decreasing per capita growth rate that results in limited growth. Each of these limited growth scenarios says something different about the process of growth or the growth relationship and each results in slightly different limited growth curves.  When the inflow of a population is exponential and the outflow is mediated by intraspecific competition for resources, $\left(\mathrm{\textit{e.g.}}, \frac{dP}{dt}=rP-\frac{r}{K}P^2\right)$, this is equivalent to saying that overall, the per capita growth rate is linear $\left(\mathrm{\textit{i.e.}}, \frac{1}{P}\frac{dP}{dt}=r\left(1-\frac{P}{K}\right)\right)$. However, growth may be limited by other processes which may result in a non-linear decrease of per capita growth rate.  In this case, one might fit many types of mathematical models to the data set and select the model that minimizes error and avoids over-fitting, for example by using a measure such as the Akaike Information Criterion (often referred to as AIC; for a review of model fitting, see \cite{Ledder16}). The model that has the lowest AIC may tell you something about which processes may be driving the population.

In some cases, the fact that a model is not matching the outcomes observed in reality can also be important. Such is the case with the primary use of the HWE model discussed above.  This model predicts the distribution of offspring genotypes in a population given a list of assumptions. If those assumptions are true, and if in the parent generation the probability or frequency of allele A occurs is $p$ and the probability that allele a occurs is $q$, then the next generation will have the following distribution of genotypes: $\mathrm{P(AA)}=p^2$, $\mathrm{P(Aa)}=2pq$, $\mathrm{P(aa)}=q^2$.  If the observed distribution in the offspring does not match, this is useful, because then we know that one of the assumptions of the HWE model have been violated.  HWE is considered a classic use-case of a null model in biology \cite{Svoboda13}.

\section{USING THE FRAMEWORK TO IMPROVE STUDENT LEARNING}

We have synthesized a framework around modeling with the view that a unified framework allows us to be more purposeful practitioners around the teaching of modeling. While traversing the modeling pathways may be intuitive for expert practitioners, some studies have found that students have difficulty translating between model representations \cite{Ainsworth99}, although the facility of translation did vary depending upon the specific relations selected \cite{Ainsworth96}.  Yet, in keeping with the recommendation from Understanding Science to be explicit \cite{Understandingsciencea}, one study showed that simply telling students the purpose of the multiple representations can have a positive impact on learning \cite{Schwonke09}.  In particular, Schwonke \textit{et al}.~conducted two studies, collecting gaze data from students viewing multiple representations of the same problem.  From the initial study, many students did not understand why they should have different representations and why they should transfer between them.  To improve the transitions between representations, one group of students in a follow-up study (reported in the same paper \cite{Schwonke09}) received additional instructions explaining the ``bridge between [the] problem texts and equations,'' while a control group did not have these instructions.  Schwonke \textit{et al}.~concluded that an explanation of the different representations improved the learning outcomes for both low- and high-prior knowledge students, but in different ways.  The low-prior knowledge group seemed to ``transfer knowledge more easily between representations,'' while the high-prior knowledge students benefited because they paid more attention to the different representations.  Our conclusion from this prior work is that while it is important for the instructor to facilitate use of multiple model representations and modeling activities, it is most valuable when the reasons behind it are made clear to the students, \textit{i.e.}, make it explicit, reflect and connect, and provide context \cite{Ainsworth96,Ainsworth99,Schwonke09,Understandingsciencea}.
%
%
%
%
%
%
``Understanding Science,'' a website that lays out a framework for the process of science, suggests the following three main actions for bringing the process of science into the classroom:  (1) make it explicit, (2) help the students reflect upon it, and (3) give it context, again and again \cite{Understandingsciencea}.  Because modeling is theoretical science, we suggest the same actions to bring our modeling framework into the classroom to strengthen both students' understanding of, and abilities in, modeling.  We list these actions again here with specific modeling framework examples:

{\bf Be explicit.}  Tell your students what you are doing, and why, in the context of modeling. Teach students the modeling framework, including definitions of models and modeling, the different model representations (boxes), and the many modeling activities (arrows).  Be clear that biologists are already engaged in modeling but may not realize it or use the same language or approach as mathematicians or statisticians. There are many first steps with which one can begin modeling, and it does not have to begin the same way that was outlined in our example modeling pathways.  Teaching any arrow can be a first step in teaching modeling as long as one is explicit in connecting it to the modeling process.

{\bf Help them reflect (and connect).}  Have students reflect upon the modeling process by assigning metacognitive exercises, such as a one-minute paper.  Work with instructors from different disciplines to help students make connections between classes to solidify their understanding of modeling.  Using a common framework can help us engage in conversations with colleagues from other disciplines, and thus bring the connections to our students.

{\bf Give it context, again and again.}  Ground the modeling investigation in the biological problem, using the {\bf Experiential} representation.  Use the framework to acknowledge and clarify the various approaches that each discipline takes to solving the same scientific problem. Show that different practitioners have different paths through the boxes that are equally valid.  Give examples of different paths in the same context.

For information and ideas around using our framework for teaching, see the following collection at \url{qubeshub.org} (\url{https://qubeshub.org/primus-ruleoffive}) \cite{QUBES15}.

\section{CONCLUDING REMARKS}
There have been numerous calls for mathematicians to work more closely with members of partner disciplines, such as biology, to improve student learning and retention in STEM (\textit{e.g., The Mathematical Sciences in 2025} and \textit{A Common Vision}, \cite{NRC13, CommonVision}).  Analogous to this, there have been numerous calls for biologists to incorporate more mathematical modeling into their curriculum to better prepare students for future careers in research and the health professions (\textit{e.g., Vision and Change}, \textit{Bio2010}, \textit{Scientific Foundations for Future Physicians}, \cite{Visionandchange,SFFP,Bio2010}).  Despite this widespread agreement on what needs to happen, implementing the recommendations has been slow.  One barrier that must be overcome to initiate true transdisciplinary conversations and collaborations on improving the modeling curriculum is to agree upon a common definition of what is meant by models and modeling between biologists, mathematicians, and statisticians.  Therefore, we have proposed adapting the ``rule of five'', which has previously informed calculus reform efforts, to describe a framework for modeling that can bring all the disciplines together. This framework defines five types of model representations (\textbf{Experiential}, \textbf{Verbal}, \textbf{Numerical}, \textbf{Visual}, and \textbf{Symbolic}) and modeling activities that provide flexible routes through a modeling investigation by students.  We give examples about what the implementation of this framework may look like in the classroom, along with the associated benefits to student learning \cite{Simundza06}.  In particular, an advantage for students working with \textbf{Experiential} representations is that they provide an entryway for understanding models that do not require the abstraction of \textbf{Symbolic} or \textbf{Visual} representations, the integration of detail needed to understand a \textbf{Numerical} representation, or the conceptual knowledge and reading comprehension needed to understand a \textbf{Verbal} representation.  Finally, we share resources that we think will be helpful for others to use.

We also hope this framework will unify practitioners coming from different parent fields (mathematics, statistics, biology, and others) and allow them to find the similarities and differences in their approaches to modeling, leading to more productive interdisciplinary conversations.  Having initiated such conversations ourselves, we have some advice based on our own experiences:  

{\bf Have a shared goal.} The goal could be simply to improve student learning in a disciplinary course, or it could be to answer a research question of mutual interest.  One of us found her way to teaching modeling through first forming a research collaboration. 
 
{\bf Start small, but have a concrete deliverable.} It is overwhelming to try to revamp a whole course; a lot can be accomplished by modifying a single lesson plan, upon which further course modifications can be built.  On the research side, aim for an internal grant proposal, a conference poster or co-advising a thesis (Deadlines help!). A few of us have been involved in using this framework to modify one lesson plan on Hardy-Weinberg Equilibrium. What lesson plan would you change first?

{\bf Be willing to be both a student and a teacher.}  Listen with respect and be open to another perspective.  We are highly trained in our own disciplines and used to being in the position of expert; learning the language, culture, and foundational knowledge of another discipline requires leaving our comfort zones, which, by definition, is uncomfortable.  Being willing to be uncomfortable requires the courage to be in the vulnerable position of learner versus expert. It can be humbling to be a student again, but is also an opportunity to remind us of what our own students experience.

{\bf Be explicit about language.} \textit{Model} is not the only word that has different meanings to different disciplines. For us, defining our language meant writing this paper. It is often necessary to clarify meanings to gain insight. One tactic is to include one or more students in on conversations. We often naturally change our language and our assumptions about prior knowledge and context to accommodate students, and this change should also benefit communication between new interdisciplinary collaborators. In addition, students benefit from observing and participating in interdisciplinary conversations. 
  
{\bf Be in it for the long term.}  Interdisciplinary relationships take time, and persistence will pay off.  If this was easy, we would not have the multiplicity of reports encouraging us to do more.

A positive outcome from these interdisciplinary conversations will be making the connections explicit to students in different disciplinary courses, reinforcing the concepts for students, and empowering them to apply knowledge from one domain to another, making them informed citizens for the 21\textsuperscript{st} century.  In their professional futures, students will not encounter textbook questions with multiple choice answers.  Instead, they will hear a wildlife biologist discuss the rate of population growth, they will see a graph in a paper that they are reading, they will monitor a population, or they will use software to run management scenario planning.  A student that can move between these representations to help solve problems is one with a superior preparation for the profession.  We have a trained disciplinary identity, but we are all modelers. We can work together to help our students be modelers, too.

\section*{ACKNOWLEDGEMENTS}
This work was conducted as part of ``Unpacking the Black Box: Teaching Quantitative Biology'' Working Group at the National Institute for Mathematical and Biological Synthesis (NIMBioS), sponsored by the National Science Foundation through NSF Award \#DBI-1300426, with additional support from The University of Tennessee, Knoxville.  Co-authors Glenn Ledder and Richard C. Schugart were supported as NIMBioS Sabbatical Fellows.  We would like to thank the other members of the working group for useful discussions: Melissa L. Aikens, Joseph T. Dauer, Samuel S. Donovan, Ben G. Fitzpatrick (also a Sabbatical Fellow at NIMBioS), Gregory D. Goins, Kristin P. Jenkins, John R. Jungck, Robert L. Mayes, and Edward F. (Joe) Redish. We would also like to thank Deborah Hughes-Hallett for her insights.  Organization of this co-authorship and working group between meetings would not have been possible without QUBES Hub supported by NSF Awards \#DBI 1346584, DUE 1446269, DUE 1446258, and DUE 1446284 \cite{QUBES15,QUBES16} and our QUBES communications student, Timothy Beaulieu.

%
%
%
%

\newpage
%
%
%
%
\section*{BIOGRAPHICAL SKETCHES}

\vspace*{.3 true cm} \noindent C.D.~Eaton earned her B.A.~in Mathematics with a minor in Zoology and her M.A.~in Interdisciplinary Mathematics from the University of Maine.  She then earned her Ph.D.~in Mathematical Ecology and Evolutionary Theory from the University of Tennessee. Dr. Eaton's primary research interests are interdisciplinary mathematics and education. She is a principal investigator on the QUBES project for Quantitative Undergraduate Biology Education and Synthesis. 

\vspace*{.3 true cm} \noindent H.c.~Highlander earned a B.A.~in Mathematics with minors in Computer Science and Music from Wesleyan College, and an M.S.~and Ph.D.~in Mathematics with an emphasis in Mathematical Biology from Vanderbilt University. She then worked as a postdoctoral fellow at the Institute for Mathematics and its Applications at the University of Minnesota.  Her primary research interests include mathematical modeling of cellular and molecular processes and infectious diseases, sensitivity analysis, and mathematical biology education.

\vspace*{.3 true cm} \noindent K.D.~Dahlquist earned a B.A.~in Biology from Pomona College and a Ph.D. in Molecular, Cellular, and Developmental Biology from the University of California, Santa Cruz. She then worked as a postdoctoral fellow at the Gladstone Institute of Cardiovascular Disease at the University of California, San Francisco.  She studies the systems biology of gene regulatory networks using genomics, mathematical, and computational biology approaches, bringing this research into the classroom.

\vspace*{.3 true cm} \noindent G.~Ledder earned a B.S.~in Ceramic Engineering from Iowa State University and a Ph.D.~in Applied Mathematics from Rensselaer Polytechnic Institute.  His research interests span a variety of theoretical science areas including groundwater flow, plant physiology, and dynamic energy budget theory.  He was the editor of a PRIMUS special issue and MAA Notes volume on mathematics pedagogy for biology students.  He is also the author of a mathematics textbook for biology students and is a member of the QUBES Advisory Board.

\vspace*{.3 true cm} \noindent M.D.~LaMar earned a B.S.~in Mathematics from the University of Texas at San Antonio and a Ph.D.~in Mathematics from the University of Texas at Austin.  His research interests include mathematical biology, specifically focusing on computational neuroscience and ecological modeling, as well as quantitative biology education.  He is also a principal investigator on the QUBES project. 

\vspace*{.3 true cm} \noindent R.C.~Schugart earned a B.A.~in Mathematics from the State University of New York College at Geneseo, and an M.S.~and Ph.D.~in Applied Mathematics from North Carolina State University.  He was a post-doctoral fellow at the Mathematical Biosciences Institute at The Ohio State University.  His research area is mathematical biology, with particular interests on modeling in wound healing and engaging undergraduate students in research.  He is also an editor for SMB Digest, a weekly moderated electronic mail containing items of current interest to the math-biology community for the Society for Mathematical Biology.


\begin{thebibliography}{0}
\bibliographystyle{plain}
\bibitem{Ainsworth96}Ainsworth S., D. Wood, and P. Bibby 1996. {\em Coordinating multiple representations in computer based learning environments}. Proceedings of the European Conference on Artificial intelligence in education. 336-342.

\bibitem{Ainsworth99}Ainsworth S. 1999. {\em The function of multiple representations}.  Computers in Education. 336-342.

\bibitem{Visionandchange} American Association for the Advancement in Science (AAS). 2010. {\em Transforming undergraduate education in biology: Mobilizing the community for change}. http://visionandchange.org. Accessed 23 May 2016.

\bibitem{AMATYC} American Mathematical Association of Two-Year Colleges (AMATYC). 2014. {\em Opening doors through mathematics. Position on professional development}. \url{http://www.amatyc.org/?page=PositionProfDev}. Accessed 26 May 2016.

\bibitem{ASA14}American Statistical Association (ASA). 2014. {\em Discovery with data: Leveraging Statistics and computer science to transform science and society}. Alexandria, VA. \url{http://www.amstat.org/policy/pdfs/BigDataStatisticsJune2014.pdf}. Accessed 31 May 2016.

\bibitem{Bender78} Bender E.A. 1978. {\em An introduction to mathematical modeling}. Wiley-Interscience. Mineola, New York. 1-46.

\bibitem{Bransford99} Bransford J.D., Brown A.L., Cocking R.R.,  Donovan M.S., and  Pellegrino J.W. (Eds.). 1999. {\em How people learn brain, mind, experience, and school}. Commission on Behavioral and Social Sciences and Education National Research Council. National  Academy Press. Washington, D.C. 1-363.

\bibitem{Brauer01} Brauer F., Castillo-Chavez C., and Bles D. 2001. {\em Mathematical models in population biology and epidemiology}. vol. 40. {em\ Springer Science Buissness Media}. LCC. Springer-Verlag, New York.

\bibitem{Carson14} Carson E.C. and Cobelli C. 2014. {\em Modeling methodology for physiology and medicine}. 2nd edition. Elsevier. Waltham, MA. 1-562. 

\bibitem{Cast11} Center for Advanced Spatial Technology (CAST). 2011. {\em Universal design for learning guidlines}. Version 2.0. Wakefeild, MA.

\bibitem{Cox95} Cox R. and Brna P. 1995. {\em Supporting the use of external representations in problem solving: The need for flexiable learning environments.} {\em Journal of Artificial Intelligence in Education}.6: 239-302.

\bibitem{Dick94} Dick T. and Patton C. 1994. {\em Calculus of a single Variable.} Boston, MA. PWS Publishing Company.

\bibitem{Douglas86} Douglas R.G. 1986. {\em Toward a lean and lively calculus: Report of the conference/workshop to develop curriculum and teaching methods for calculus at the college level. Journal of College Mathematics}. 18(5): 439-442.

\bibitem{Donovana15} Donovan S., Eaton C.D., Gower S.T., {\em et al}. 2015. {\em QUBES: A community focused on supporting teaching and learning in quantitative biology. Letter in Biomathematics}.2(1): 46-55. 

\bibitem{Epstein} Epstein, Joshua M. 2008. {\em Why Model?}. Journal of Artificial Societies and Social Simulation. 11(4)12. \url{http://jasss.soc.surrey.ac.uk/11/4/12.html}.  Accessed 21 March 2017.

\bibitem{Ganter04} Ganter S. and Barker W. (EDs.) 2004. {\em The curriculum foundations project: Voices of the partner disciplines}. Washington, DC: Mathematics Association of America.

\bibitem{Garfunkel16} Garfunkel S. and Montgomery M. (EDs.) 2016. {\em Guideline for assessment and Instruction in mathematical modeling education}. Bedford, MA. Consortium for Mathematics and Its Applications. 

\bibitem{Gleason92} Gleason A.M. and Hughes-Hallet D. 1992. {\em The Calculus Consortium Based at Harvard University}. In \emph{Focus on Calculus: A Newsletter for the Calculus Consortium at Harvard University}. Issue 1. Harvard University.

\bibitem{Hiebert86} Hiebert J. 1986. {\em Conceptual and procedural knowledge: The case of mathematics.} Mahwah, NJ. Erlbaum Associates. 

\bibitem{Hiebert92} Heibert J. and Carpenter T.P. 1992. {\em Learning and teaching with understanding}. In Grouws D.A. (ED.) {\em Handbook oof Research on Mathematics Teaching and Learning.} New York, NY. Macmillian. 

\bibitem{Holling59a} Holling C.S. 1959a. {\em The components of predation as revealed by a study of small-mammal predation of the European sawfly. The Canadain Entomologist}. 91: 293-320.

\bibitem{Holling59b} Holling C.S. 1959b. {\em Some characteristics of simple types of predation and parasitism. The Canadian Entomologist}. 91: 293-320. 

\bibitem{SFFP} Howard Hughes Medical Institute Committee (HHMI). 2009. {\em Scientific Foundations for Future Physicians}. Report of the American Association of Medical Colleges. University of Delaware. 1-46.

\bibitem{Hughes-Hallett06} Hughes-Hallett D. 2006. {\em What Have We Learned from Calculus Reform? The Road to Conceptual Understanding}. in Hastings N. (Ed.) {\em Rethinking the Courses Below Calculus}. Mathematics Association of America.

\bibitem{Hughes-Hallett16} Hughes-Hallett D. 2016. Personal Interview. May 31, 2016.

\bibitem{Hughes-Hallett95} Huges-Hallett D., Gleason A.M., McCallum W.G., {\em et al}. 1995. {\em Calculus: Single and multivariate}. New York, NY. Wiley and Sons.

\bibitem{Hughes-Hallett98} Hughes-Hallett D., Gleason A.M., McCallum W.G., {\em et al}. 1998. {\em Calculus: Single and multivariate}. 2nd edition. New York, NY. Wiley and Sons. 

\bibitem{Hughes-Hallett05} Hughes-Hallett D., Gleason A.M., McCallum W.G., {\em et al}. 2005. {\em Calculus: Single and multivariate}. 4th edition. New York, NY. Wiley and Sons.

\bibitem{Hughes-Hallett13} Hughes-Hallett D., Gleason A.M., McCallum W.G., {\em et al}. 2013. {\em Calculus: Single and multivariate}. 6th edition. New York, NY. Wiley and Sons.

\bibitem{Jungck10} Jungck J.R., Gaff H., and Weisstein A.E. 2010. {\em Mathematical manipulative models: In defense of ``Beanbag Biology''}.  CBE-Life Sciences Education. 9: 201-211. 

\bibitem{Kolb05} Kolb A. and Kolb D.A. 2005. {\em Learning styles and learning spaces: Enhancing experiential learning in higher education. Academy of Management Learning and Education}. 4(2): 193-212.

\bibitem{Kolb15} Kolb D.A. 2015. {\em Experiential learning: Experience as the source of learning and development}. 2nd edition. Prentice Hall. Englewood Cliffs, NJ. 1-377.

\bibitem{Kolb01} Kolb D.A., Boyatiz R.E., Mainemelis C. (Eds.) 
Sternberg S. and Zhang L. 2001. {\em Experiential learning theory: previous research and new directions}. Lawrence Erlbaum Associates. New York, NY. 197-226. 

\bibitem{Ledder05} Ledder G. 2005. {\em Differential equations: A modeling approach}. McGraw-Hill Higher Education. University of Nebraska. 1-136.

\bibitem{Bugbox} Ledder G. 2007. {\em BUGBOX-predator}. \url{https://www.math.unl.edu/~gledder1/BUGBOX/}. Accessed 26 May 2016.

\bibitem{Ledder16} Ledder G. 2016. {\em Emperical modeling: Choosing models and fitting them to data.} {\em The College Mathematics Journal.} 47: 1-11.

\bibitem{CRAFTY} Mathematical Association of America (MAA). 2004. {\em Curriculum renewal across the first two years (CRAFTY)}. \url{http://www.maa.org/programs/faculty-and-departments/curriculum-department-guidelines-recommendations/crafty}. Accessed 23 May 16.

\bibitem{Bio2010} National Academy of Science (NAS). 2003. {\em BIO 2010: Transforming undergraduate education for the future research biologists}. Washington, DC. The National Academic Press. 1-192.

\bibitem{NCTM89} National Council of Teachers of Mathematics (NCTM). 1989. {\em Curriculum and evaluation standards for school mathematics.} Reston, VA. The National Council of Teachers of Mathematics. 

\bibitem{NRC13}  National Research Council (NRC). 2013.{\em The mathematical sciences in 2025}. Washington, DC. The National Academies Press. 1-222.

\bibitem{Nandi13} Nandi B. and Chakravarti A. 2013. {\em An introduction to mathematical modeling}. Jeevanu Times. 13(1):16-18.

\bibitem{Odenbaugh}	Odenbaugh J. 2005. {\em Idealized, Inaccurate but Successful: A Pragmatic Approach to Evaluating Models in Theoretical Ecology}. Biol Philos. Kluwer Academic Publishers. 20(2-3):231–55.

\bibitem{Ostebee94} Ostebee, A. and Zorn, P. 1994. {\em Calculus from graphical, numerical, and symbolic points of view.} Orlando, FL. Saunders College Publishing.  

\bibitem{PCAST12} Presidents Council of Advisers on Science and Technology (PCAST). 2012. {\em Engage to excel: Producing one million additional college graduates with degrees in science, technology, engineering, and mathematics}. 1-130.

\bibitem{QUBES15} Quantitative Undergraduate Biology Education and Synthesis (QUBES). 2015. {\em NIMBioS working group: Unpacking the black box.} \url{https://qubeshub.org/groups/nimbios_wg_teachingquantbio/overview}. Accessed 26 May 2016. 

\bibitem{QUBES16} Quantitative Undergraduate Biology Education and Synthesis (QUBES). 2016.  {\em QUBES: The Power of Biology$\times$Math$\times$Community}. \url{https://qubeshub.org}. Accessed 26 May 2016.

\bibitem{Raven16} Raven P. and Johnson G. 2016. {\em Biology}. McGraw-Hill Education. 1408.

\bibitem{Redish15} Redish F.E. and Kuo E. 2015. {\em Language of Physics, Langauge of Math: Disciplinary culture and dynamic epistemology. Science and education}. 24(5): 561-590.

\bibitem{Saxe1873} Saxe J.G. 1873. {\em The Poems}. Highgate edition. Volume 1.  James R. Osgood and Company. Boston, MA. 1-485. 

\bibitem{CommonVision} Saxe K. and Braddy L. 2015. {\em A common vision for undergraduate mathematical sciences programs in 2025}. Washington, DC. Mathematics Association of America. 

\bibitem{Schwonke09} Schwonke R., Berthold K. and Renkl A. 2009. {\em How multiple external representations are used and how they can be made more useful. Applied Cognitive Psychology} 23: 1227-1243. 

\bibitem{Simundza06} Simundza G. 2006. {\em The fifth rule: Direct experience of mathematics, a fresh start for collegiate mathematics}. Mathematical Association of America. 320-328. 

\bibitem{SIAM12} Society for Industrial and Applied Mathematics (SIAM). 2012. {\em Modeling across the curriculum I.}. Philidelphia, PA. 

\bibitem{SIAM14} Society for Industrial and Applied Mathematics (SIAM). 2014. {\em Modeling across the curriculum II.}. Philidelphia, PA.

\bibitem{Sorgo10} Sorgo A. and Jungck J. 2010. {\em Connecting biology and mathematics: First Prepare the teachers. CBE Life Sciences Education}. 9(3):196-200. 

\bibitem{Steele16} Steele J. and Demand Media. 2016. {\em Steady state vs. equilibrium in Biology}. Seattle PI.  \url{http://education.seatllepi.com/steady-state-vs-equilibruim-biology-6085.html}. Accessed 19 May 2016.

\bibitem{Svoboda13} Svoboda J. and Passmore C. 2013. {\em The strategies of modeling in biology education. Science and Education.} 22(1): 119-142.

\bibitem{Tabachneck94} Tabachneck H.J.M., Koedinger K.R. and Nathan M.J. 1994. {\em Towards a theoretical account of strategy use and sense making in mathematical problem solving}. Proceedings of the 16th annual conference of the cognitive science society. Hillsdale, NJ. Erlbaum. 836-841.

\bibitem{Tufte}  Tufte E.R. 2001. {\em The Visual Display of Quantitative Information}.  2nd edition. Graphics Press.

\bibitem{Understandingsciencea} Understanding Science. 2016a. {\em University of California museum of paleontology}. \url{http://www.understandingscience.org}. Accessed 19 May 2016.

\bibitem{UnderstandingScienceb} Understanding Science. 2016b. {\em The real process of science}. \url{http://undsci.berkeley.edu/article/howscienceworks_02}. Accessed 19 May 2016.

\bibitem{Utter96} Utter F.W. 1996. {\em Relationship among AP calculus teachers' pedagogical content beliefs, classroom practice and their student’s acheivement}. Dissertation, Department of Mathematics. Oregon State University. 

\bibitem{Waldron16} Waldron I. 2016. {\em Population growth-exponential and logistic models vs. complex reality}. University of Pennsylvannia. Philadelphia, PA. \url{http://serendip.brynmawr.edu/exchange/bioactivities/pop}. Accessed 19 May 2016.

\bibitem{Watson53} Watson J.D. and Crick F.H. 1953. {\em Molecular structure of nucleic acids}. {\em Nature}. 171(4356): 737-738.

\bibitem{Wilensky99} Wilensky U. 1999. NetLogo. {\em Center for connected learning and computer-based modeling}. Northwestern University. Evanston, IL. \url{http://ccl.northwestern.edu/netlogo/}. Accessed 19 May 2016.

\bibitem{Windham08} Windham D.M., 2008. {\em Faculty perceptions of a calculus reform experiment at a research university: A historical qualitative analysis}. Dissertation, Department of Middle and secondary education, Florida State University. ProQuest LLC. \url{http://diginole.lib.fsu.edu/islandora/object/fsu%3A169163}. Accessed 19 May 2016. 

\end{thebibliography}
\end{document}